\documentclass[12pt]{article}

\usepackage{amsfonts}
\usepackage{amsmath}
\usepackage{amsthm} 
\usepackage{epsfig}
\usepackage{graphics}
\usepackage{pslatex}

\usepackage{latexsym}
\usepackage{amssymb,exscale}
\usepackage[all]{xy}
\usepackage{hyperref}
\usepackage{verbatim} 

\theoremstyle{plain} 
    \newtheorem{theorem}{Theorem}
    \newtheorem{lemma}[theorem]{Lemma}
    \newtheorem{proposition}[theorem]{Proposition}
    \newtheorem{corollary}[theorem]{Corollary}

\theoremstyle{definition} 
    \newtheorem{definition}[theorem]{Definition}

    \newtheorem{remark}[theorem]{Remark}

\def\eps{\epsilon}

\def\ra{\rightarrow}

\def\C{\mathbb{C}}

\def\R{\mathbb{R}}

\def\e{{\varepsilon}}

\def\z{{\bf z}}

\def\<{\langle}
\def\>{\rangle}

\newcommand{\E}{\mbox{\bf E}}

\newcommand\tr{{\mbox{\rm tr}}}

\newcommand\floor[1]{\left\lfloor #1 \right\rfloor}
\newcommand\norm[1]{\left\| #1 \right\|}
\newcommand\ceiling[1]{\lceil #1 \rceil}
\newcommand\mnote[1]{} 
\newcommand\be{\begin{equation*}}

\newcommand\ee{\end{equation*}}

\newcommand\ben{\begin{equation}}
\newcommand\een{\end{equation}}
\newcommand\bes{\begin{eqnarray*}}
\newcommand\ees{\end{eqnarray*}}

\newcommand{\sm}{{\raise0.3ex\hbox{$\scriptstyle \setminus$}}}

\def\eps{\epsilon}

\renewcommand{\phi}{\varphi}

\newcommand{\tbn}{T_{b,N}}

\def\CHI{\mathchoice%
{\raise2pt\hbox{$\chi$}}%
{\raise2pt\hbox{$\chi$}}%
{\raise1.3pt\hbox{$\scriptstyle\chi$}}%
{\raise0.8pt\hbox{$\scriptscriptstyle\chi$}}}
\def\smalloplus{\raise1pt\hbox{$\,\scriptstyle \oplus\;$}}

\author{Alice Guionnet\thanks{
UMPA, CNRS  UMR 5669, ENS Lyon, 46 all\'ee d'Italie,
69007 Lyon, France. {aguionne@umpa.ens-lyon.fr}.
This work was partially supported by the ANR project
ANR-08-BLAN-0311-01.}{ } \;
Philip Wood \;
 Ofer Zeitouni\thanks{School of Mathematics,
University of Minnesota and Faculty of Mathematics,
Weizmann
Institute, POB 26, Rehovot 76100, Israel. {zeitouni@math.umn.edu}. 
The work of this author was partially
supported by NSF grant DMS-0804133 and by a grant from the Israel 
Science Foundation.} }
\title{Convergence of the spectral measure of  non normal matrices}

\begin{document}
 \bibliographystyle{abbrv}
\date{ {\small
October 11, 2011}}
\maketitle
\abstract{We discuss regularization by noise of the spectrum of large random
non-Normal matrices. Under suitable conditions, we show that the regularization
of a sequence of matrices that converges in $*$-moments to a regular
element $a$, 
by the addition of a polynomially vanishing Gaussian Ginibre matrix, forces 
the empirical measure of eigenvalues
to converge to the Brown measure of $a$.}
\section{Introduction}
Consider a sequence $A_N$ of $N\times N$ matrices, of uniformly bounded
operator norm, and assume that
$A_N$ converges in $*$-moments toward an element $a$ in a $W^*$ probability
space $({\cal A}, \|\cdot\|, *, \phi)$, that is, for any non-commutative 
polynomial $P$,
$$ \frac1N \mbox{\rm tr} P(A_N,A_N^*)\to_{N\to\infty} 
\phi(P(a,a^*))\,.$$
We assume throughout that the tracial state $\phi$ is faithful; this
does not represent a loss of generality.
If $A_N$ is a sequence of Hermitian matrices, this is enough in order to
conclude that the empirical measure of eigenvalues of $A_N$, that is
the measure
$$ L_N^A:=\frac1N \sum_{i=1}^N \delta_{\lambda_i(A_N)},$$
where $\lambda_i(A_N), i=1\ldots N$ are the eigenvalues of $A_N$,
converges weakly to a
 limiting measure $\mu_a$,
  the spectral measure of $a$, 
 supported on a compact subset of
 $\R$. (See \cite[Corollary 5.2.16, Lemma 5.2.19]{AGZ} for this standard
 result and further background.) Significantly, in the Hermitian case,
 this convergence is stable under
 small bounded perturbations: with $B_N=A_N+E_N$ and $\|E_N\|<\epsilon$,
 any subsequential limit of $L_N^B$ will belong to 
 $B_L(\mu_a,\delta(\epsilon))$, with $\delta(\epsilon)\to_{\epsilon\to 0} 0$
 and $B_L(\nu_a,r)$ is the ball (in say, the L\'{e}vy metric)
 centered at $\nu_a$ and of radius $r$.

 Both these statements fail when $A_n$ is not self adjoint. For a standard
 example (described in \cite{Sniady}),
 consider the nilpotent matrix
$$
T_N=\left(\begin{array}{lllll}
  0&1&0&\ldots&0\\
  0&0&1&0&\ldots\\
\ldots&\ldots&\ldots&\ldots&\ldots\\
0&\ldots&\ldots&0&1\\
0&\ldots&\ldots&\ldots&0
\end{array}
\right)\,.$$
Obviously, $L_N^T=\delta_0$, while a simple computation
reveals that $T_N$ converges in $*$-moments to
a Unitary Haar element of ${\cal A}$, 
that is
\begin{equation}
  \label{eq-tn}
  \frac1N \mbox{\rm tr}(T_N^{\alpha_1} (T_N^*)^{\beta_1}
\ldots T_N^{\alpha_k} (T_N^*)^{\beta_k})
\to_{N\to\infty} \left\{\begin{array}{ll}
  1, & \mbox{\rm if} \sum_{i=1}^k\alpha_i=\sum_{i=1}^k\beta_i,\\
  0, & \mbox{\rm otherwise}.
\end{array}
\right.
\end{equation}
Further,  adding to $T_N$ the matrix whose entries are all $0$ except for the 
bottom left, which is taken as $\epsilon$, changes the empirical measure 
of eigenvalues drastically - as we will see below, as $N$ increases, 
the empirical measure converges to the uniform measure on the unit circle in
the complex plane.

Our goal in this note is to explore this phenomenun in the context of 
small random perturbations of matrices.
We recall some notions. 
For $a\in {\cal A}$, the {\it Brown measure} $\nu_a$ on $\C$ 
is the measure satisfying
$$ \log \mbox{\rm det}(z-a)=\int \log |z-z'| d\nu_a(z'),\quad
\z\in \C,$$
where $\mbox{\rm det}$ is the Fuglede-Kadison determinant; we refer to
\cite{brown,haageruplarsen} for  definitions. We have in particular that 
$$
\log \mbox{\rm det}(z-a)=\int \log x d\nu_{a}^z(x)\,\quad z\in \C\,,$$ 
where $\nu_a^z$ denotes the spectral measure of the operator $|z-a|$.
In the sense of distributions, we have
$$ \nu_a=\frac{1}{2\pi} \Delta \log \mbox{\rm det}(z-a)\,.$$
That is, for smooth compactly supported function $\psi$ on $\C$,
\begin{eqnarray*}
  \int  \psi(z) d\nu_a(z)&=&
  \frac{1}{2\pi}\int dz \ \Delta \psi(z)
  \int \log|z-z'| d\nu_a(z')\\
  &=&
  \frac{1}{2\pi}\int dz\ \Delta \psi(z)
  \int \log x d\nu_a^z (x)\,.
\end{eqnarray*}
A crucial assumption in our analysis is the following.
\begin{definition}[Regular elements]
  An element $a\in {\cal A}$ is {\em regular} if
  \begin{equation}
    \label{eq-3}
    \lim_{\epsilon\to 0} 
    \int_{\C}dz \Delta \psi(z) \int_0^\epsilon \log x d\nu_a^z(x)
=0\,,
\end{equation}
for all smooth functions $\psi$ on $\C$ with compact support.
\end{definition}
Note that regularity is a property of $a$, not merely of its Brown measure
$\nu_a$.
We next introduce the class of Gaussian perturbations we consider.
\begin{definition}[Polynomially vanishing Gaussian matrices]
A sequence of $N$-by-$N$ random Gaussian matrices is called 
{\em polynomially vanishing} if its entries $(G_N(i,j))$
are independent
centered complex Gaussian variables, and there exist $\kappa>0$, $\kappa'
\geq 1+\kappa$ so that
$$ N^{-\kappa'}\leq E|G_{ij}|^2 \leq N^{-1-\kappa}\,.$$
\end{definition}
\begin{remark}
\label{rem-1}
As will be clear below, see the beginning of the proof of
Lemma \ref{lem-5}, the Gaussian assumption only intervenes in obtaining a 
uniform lower bound on singular values of certain random matrices. 
As pointed out to us by R. Vershynin,
this uniform estimate extends to other situations, most notably to 
the  polynomial rescale of matrices whose entries are i.i.d. and possess 
a bounded density.
We do not discuss such extensions here.
\end{remark}

Our first result is a stability,
with respect to polynomially vanishing
Gaussian
perturbations, of the convergence of spectral measures for non-normal matrices.
Throughout, we denote by $\|M\|_{op}$ the operator norm 
of a matrix $M$.
\begin{theorem}\label{theo1}
  Assume that the uniformly bounded (in the operator norm) 
  sequence of $N$-by-$N$ matrices
 $A_N$ converges in $*$-moments to a regular element $a$.
Assume further that $L_N^A$ converges weakly to the Brown measure
$\nu_a$.
Let $G_N$ be a sequence of polynomially vanishing Gaussian
matrices, and set $B_N=A_N+G_N$.
Then,  $L_N^{B}\to \nu_a$ weakly, in probability.
\end{theorem}
Theorem \ref{theo1} puts rather stringent assumptions on the sequence $A_N$.
In particular, its assumptions are not satisfied by the sequence of 
nilpotent matrices $T_N$ in \eqref{eq-tn}. Our second 
result corrects this defficiency, by showing that small Gaussian 
perturbations ``regularize'' matrices that are close to
matrices satisfying the assumptions of Theorem \ref{theo1}.
\begin{theorem}\label{theo2}
Let $A_N$, $E_N$ be a sequence of 
bounded (for the operator norm) $N$-by-$N$ matrices, so that
$A_N$ converges in $*$-moments to a regular element $a$.
Assume that $\|E_N\|_{op}$ converges to zero polynomially fast
in $N$, and that $L_N^{A+E}\to \nu_a$ weakly. Let
$G_N$ be a sequence of polynomially vanishing Gaussian matrices, and
set $B_N=A_N+G_N$. Then, $L_N^B\to \nu_a$ weakly, in probability.
\end{theorem}
Theorem \ref{theo2} should be compared to earlier 
results of Sniady \cite{Sniady}, who used stochastic calculus to show
that a perturbation by an asymptotically vanishing Ginibre Gaussian
matrix regularizes arbitrary matrices. Compared with his results, we allow for
more general Gaussian perturbations (both structurally and in terms
of the variance) and also show that the Gaussian regularization can decay as
fast as wished in the polynomial scale. On the other hand, we do impose a
regularity property on the
limit $a$ as well as on the sequence of matrices for
which we assume that adding a polynomially small matrix is
enough to obtain convergence to the Brown measure.  

A  corollary of our general results is the following.
\begin{corollary}
  \label{cor-tn}
  Let $G_N$ be a sequence of polynomially vanishing Gaussian matrices and let
  $T_N$ be as in \eqref{eq-tn}. Then $L_N^{T+G}$ converges weakly, in 
  probability, 
  toward the uniform measure on the unit circle in $\C$.
\end{corollary}

\begin{figure}
\begin{center}
        \fbox{
\begin{minipage}{\textwidth}
\begin{center}
\begin{tabular}{cccc}
\hspace{.7cm}\parbox[t]{1.2in}{
{\bf (a) } \parbox[t]{.8in}{$N=50$ }\\[1pt]
\scalebox{.55}{\hspace{-1.5cm}\includegraphics{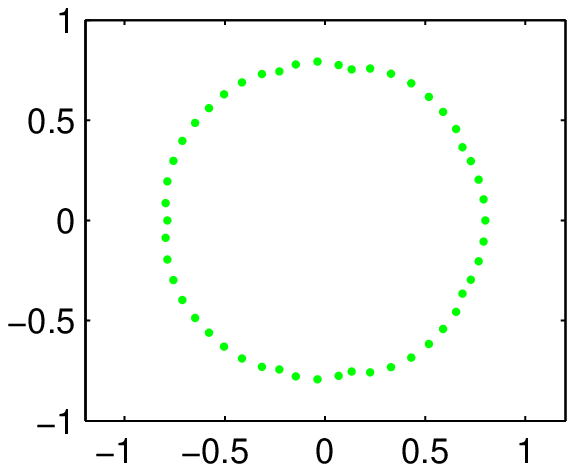}}
}
&
\parbox[t]{1.2in}{
{\bf (b) } \parbox[t]{.8in}{$N=100$}\\[1pt]
\scalebox{.55}{\hspace{-1.5cm}\includegraphics{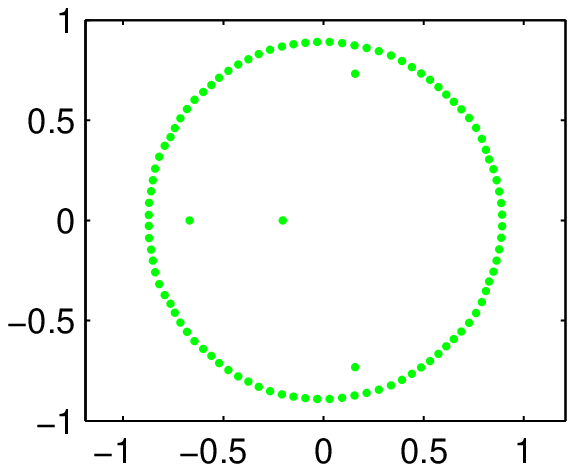}}
}
&
\parbox[t]{1.2in}{
{\bf (c) } \parbox[t]{.8in}{$N=500$}\\[1pt]
\scalebox{.55}{\hspace{-1.5cm}\includegraphics{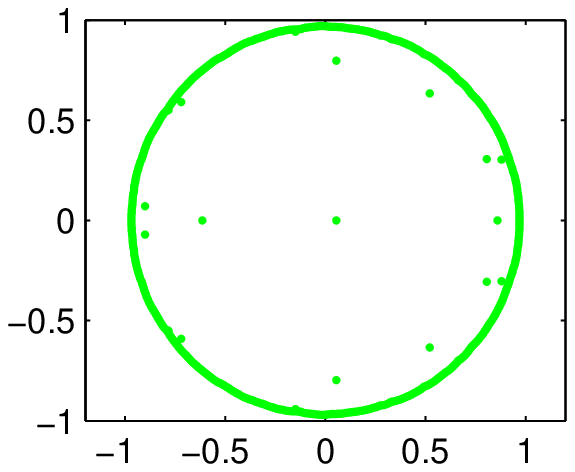}}
}
&
\parbox[t]{1.2in}{
{\bf (d) } \parbox[t]{.8in}{$N=5000$}\\[1pt]
\scalebox{.55}{\hspace{-1.5cm}\includegraphics{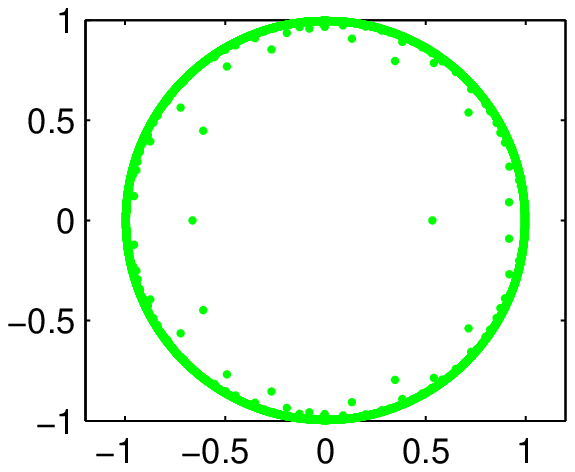}}
}
\end{tabular}
\end{center}
\vspace{-.5cm}
\caption{The eigenvalues of $T_N+N^{-3-1/2}G_N$, where
$G_N$ is iid complex Gaussian with mean 0, variance 1 entries.
}
\label{fig2}
\end{minipage}\ \ \ 
}
\end{center}
\vspace{-.8cm}
\end{figure}

\noindent
In Figure~\ref{fig2}, we give a simulation of the setup in
Corollary~\ref{cor-tn} for
various $N$.

We will now define class of matrices $T_{b,N}$ for which, if $b$ is chosen correctly, adding a polynomially vanishing Gaussian matrix $G_N$ is not sufficient to regularize $T_{b,N}+G_N$.  
Let $b$ be a positive integer, and define $T_{b,N}$ to be an $N$ by $N$ block diagonal matrix which each $b+1$ by $b+1$ block on the diagonal equal $T_{b+1}$ (as defined in \eqref{eq-tn}.  If $b+1$ does not divide $N$ evenly, a block of zeros is inserted at bottom of the diagonal.
Thus, every entry of $T_{b,N}$ is zero except for entries on the superdiagonal (the superdiagonal is the list of entries with coordinates $(i,i+1)$ for $1\le i\le N-1$), and the superdiagonal of $T_{b,N}$ is equal to
$$(\underbrace{1,1,\dots,1}_b,0,\underbrace{1,1,\dots,1}_b,0,\dots,\underbrace{1,1,\dots,1}_b,\underbrace{0,0,\dots,0}_{\le b}).$$
Recall that the spectral radius of a matrix is the maximum absolute value of the eigenvalues.  Also, we will use $\|A \|=\tr(A^*A)^{1/2}$ to denote the Hilbert-Schmidt norm.

\begin{proposition}\label{prop-lb}
   Let $b= b(N)$ be a sequence of positive integers such that $b(N)\ge log N$ for all $N$, and let $T_{b,N}$ be as defined above. Let $R_N$ be an $N$ by $N$ matrix satisfying $\| R_N\| \le g(N)$, where for all $N$ we assume that $g(N) < \frac1{3 b \sqrt N}$.  Then
$$ \rho( T_{b,N} + R_N) \le (Ng(N))^{1/b} +o(1),$$
where $\rho(M)$ denotes the spectral radius of a matrix $M$, and $o(1)$ denotes a small quantity tending to zero as $N\to \infty$.
\end{proposition}

Note that $T_{b,N}$ converges in $*$-moments to 
a Unitary Haar element of $\mathcal A$ (by a computation similar to \eqref{eq-tn}) if $b(N)/N$ goes to zero, which is a regular element.  The Brown measure of the Unitary Haar element is uniform measure on the unit circle; thus, in the case where $(Ng(N))^{1/b}<1$, Proposition~\ref{prop-lb} shows that $\tbn+R_N$ does not converge to the Brown measure for $\tbn$.

\begin{corollary}\label{cor-lb}
Let $R_N$ be an iid Gaussian matrix where each entry has mean zero and variance one.   
 Set $b=b(N)\ge \log N$ be a sequence of integers, and let $\gamma >5/2$ be a constant.
Then, with probability tending to 1 as $N\to \infty$, we have
$$\rho(\tbn + \exp(-\gamma b) R_N) \le \exp\left(-\gamma + \frac{2\log N}{b}\right) + o(1),$$
where $\rho$ denotes the spectral radius and where $o(1)$ denotes a small quantity tending to zero as $N\to\infty$.  Note in particular that the bound on the spectral radius is strictly less than $\exp(-1/2)< 1$ in the limit as $N\to\infty$, due to the assumptions on $\gamma$ and $b$.
\end{corollary}

\begin{figure}
\begin{center}
        \fbox{
\begin{minipage}{\textwidth}
\begin{center}
\begin{tabular}{cccc}
\hspace{.7cm}\parbox[t]{1.2in}{
{\bf (a) } \parbox[t]{.8in}{$N=50$ }\\[1pt]
\scalebox{.55}{\hspace{-1.5cm}\includegraphics{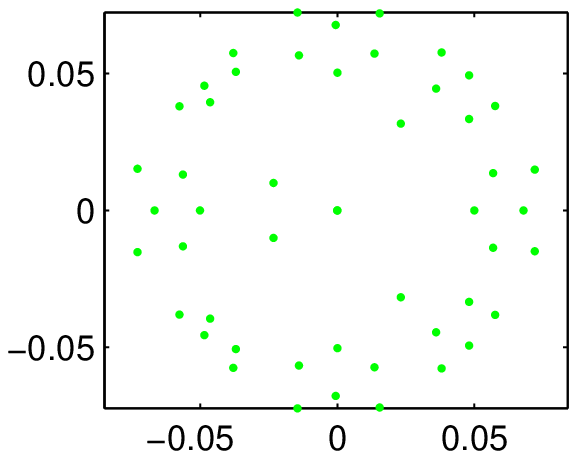}}
}
&
\parbox[t]{1.2in}{
{\bf (b) } \parbox[t]{.8in}{$N=100$}\\[1pt]
\scalebox{.55}{\hspace{-1.5cm}\includegraphics{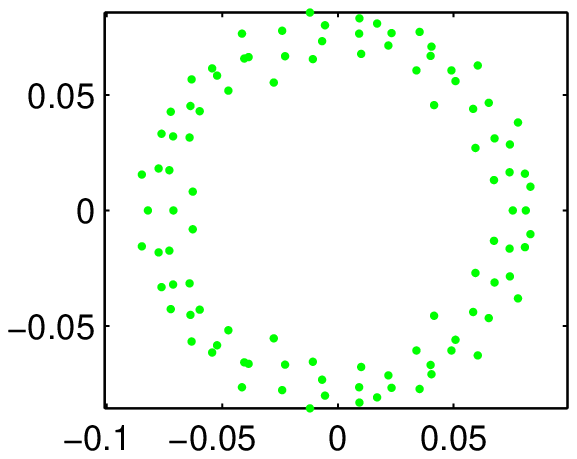}}
}
&
\parbox[t]{1.2in}{
{\bf (c) } \parbox[t]{.8in}{$N=500$}\\[1pt]
\scalebox{.55}{\hspace{-1.5cm}\includegraphics{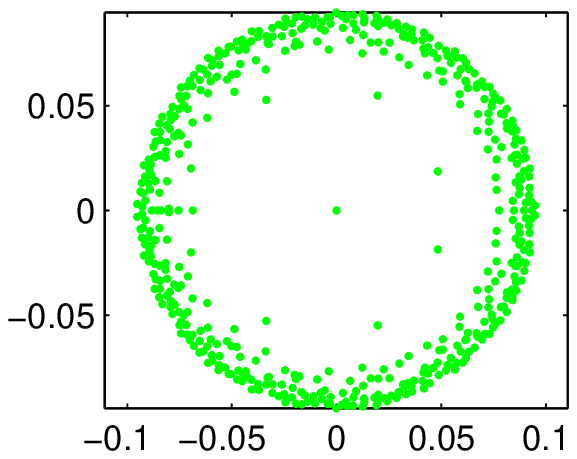}}
}
&
\parbox[t]{1.2in}{
{\bf (d) } \parbox[t]{.8in}{$N=5000$}\\[1pt]
\scalebox{.55}{\hspace{-1.5cm}\includegraphics{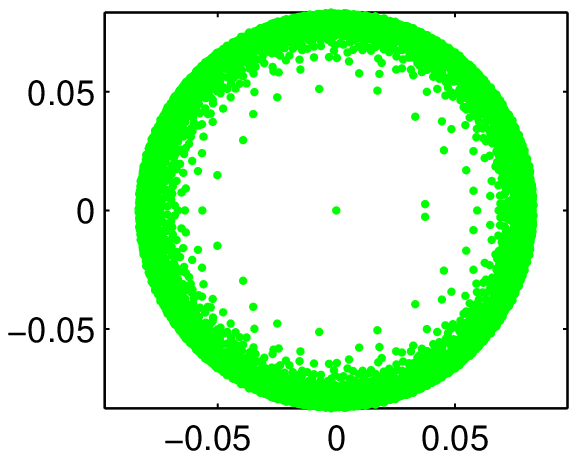}}
}
\end{tabular}
\end{center}
\vspace{-.5cm}
\caption{The eigenvaules of $T_{\log N,N}+N^{-3-1/2}G_N$, where
$G_N$ is iid complex Gaussian with mean 0, variance 1 entries.  The spectral radius
is roughly $0.07$, and the bound from Corollary~\ref{cor-lb} is
$\exp(-1)\approx 0.37$.}
\label{fig3}
\end{minipage}\ \ \ 
}
\end{center}
\vspace{-.8cm}
\end{figure}

Corollary~\ref{cor-lb} follows from Proposition~\ref{prop-lb} by noting that, with probability tending to 1, all entries in $R_N$ are at most $C\log N$ in absolute value for some constant $C$, and then checking that the hypotheses of Proposition~\ref{prop-lb} are satisfied for $g(N)=\exp(-\gamma b) CN(\log N)^{1/4}$.  There are two instances of Corollary~\ref{cor-lb} that are particularly interesting: when $b=N-1$, we see that a exponentially decaying Gaussian perturbation does not regularize $T_N=T_{N-1,N}$, and when $b=\log(N)$,
we see that polynomially decaying Gaussian perturbation does not regularize
$T_{\log N,N}$ (see Figure~\ref{fig3}).

We will prove Proposition~\ref{prop-lb} in Section~\ref{sec-proppf}.  The proof of our main results (Theorems~\ref{theo1} and \ref{theo2}) borrows from the methods of \cite{GKZ}. We introduce notation. For any $N$-by-$N$ matrix $C_N$, let
 $$\widetilde C_N=\left(
\begin{array}{cc}
0& C_N\cr
C_N^*&0\cr\end{array}
\right).$$ 
 We denote by $G_{C}$  the Cauchy-Stieltjes transform of the spectral
measure of the matrix $\widetilde C_N$, that is
$$ G_C(z)=\frac1{2N}\tr(z-\widetilde C_N)^{-1}\,, \quad z\in \C_+\,.$$
The following estimate is immediate from the definition and the resolvent
identity:
\begin{equation}
	\label{eq-comp}
	|G_C(z)-G_D(z)|\leq \frac{ \|C-D\|_{op}}
	{|\Im z|^2}\,.
\end{equation}

\section{Proof of Theorem \ref{theo1}}
We keep throughout the notation and assumptions of the theorem. The following
is a crucial simple observation.
\begin{proposition}\label{prop1}
For all complex number $\xi$, 
 and 
all $z$ so that $\Im z\ge N^{-\delta}$ with $\delta<\kappa/4$,
$$\E|\Im G_{B_N+\xi}(z)|\le \E|\Im G_{A_N+\xi}(z)|+1$$
\end{proposition}
\proof Noting that 
\begin{equation}
\label{eq-doublestar}
\E\|B_N-A_N\|_{op}^k=\E\|G_N\|_{op}^k\leq C_k N^{-\kappa k/2},
\end{equation}
the conclusion 
follows from \eqref{eq-comp} and H\"{o}lder's inequality.
\qed

We continue with the proof of Theorem \ref{theo1}. Let $\nu_{A_N}^z$ denote
the empirical measure of the eigenvalues of the matrix $\widetilde{A_N-z}$.
We have that, for smooth test functions $\psi$,
$$ \int dz \Delta \psi(z) \int \log |x| d\nu_{A_N}^z(x)=
\frac{1}{2\pi}\int \psi(z) dL_N^A(z)\,.$$
In particular, the convergence of $L_N^A$ toward $\nu_a$ implies that
$$ \E\int dz \Delta \psi(z) \int \log |x| d\nu_{A_N}^z(x)
\to \int \psi(z)d\nu_a(z)=
\int dz \Delta \psi(z) \int \log x d\nu_a^z(x)\,.
$$
On the other hand, since $x\mapsto \log x$ is bounded continuous on
compact subsets of $(0,\infty)$, it also holds that 
for
any continuous bounded function $\zeta: \R_+\mapsto \R$ compactly supported
in $(0,\infty)$,
$$ \E\int dz \Delta \psi(z) \int \zeta(x)\log x d\nu_{A_N}^z(x)
\to \int dz \Delta \psi(z)\int \zeta(x) \log x 
d\nu_a^z(x)\,.$$
Together with  the fact that $a$ is regular and that $A_N$ 
is uniformly bounded, one concludes therefore that
%
%
$$\lim_{\e\downarrow 0}\lim_{N\ra\infty}
\E \int \int_0^\e \log |x| d\nu_{A_N}^z(x)dz=0\,.$$
%
%
%
Our next goal is to show that the same applies to $B_N$. In the following,
we let $\nu_{B_N}^z$ denote the empirical measure
of the eigenvalues of $\widetilde{B_N-z}$.
\begin{lemma}
  \label{lem-5}
$$\lim_{\e\downarrow 0}\lim_{N\ra\infty}
\int \E[ \int_0^\e \log |x|^{-1} d\nu_{B_N}^z(x) ]dz=0$$
\end{lemma}
Because 
$\E\|B_N-A_N\|_{op}^k\to 0$ for any $k>0$, we have 
for any fixed smooth $w$ compactly supported in
$(0,\infty)$ that
$$ \E|\int dz \Delta \psi(z) \int w(x)\log x d\nu_{A_N}^z(x)
- \int dz \Delta \psi(z) \int w(x)\log x d\nu_{B_N}^z(x)|
\to_{N\to\infty} 0\,,$$
Theorem \ref{theo1} follows at once from Lemma \ref{lem-5}.

\noindent
{\bf Proof of lemma \ref{lem-5}:}
Note first 
that  by \cite[Theorem 3.3]{SST} (or its generalization in \cite[Proposition
16]{GKZ}
to the complex case), there exists 
a constant $C$ so that  for any $z$, the smallest singular value $\sigma^z_N$
of $B_N+zI$ satisfies 
$$P(\sigma^z_N\le x)\le C\left( N^{\frac{1}{2}+\kappa'}
x\right)^\beta$$
with $\beta=1$ or $2$ according whether we
are in the real or the complex case. 
Therefore, for any $\zeta>0$, uniformly in $z$
\begin{eqnarray*}
\E[ \int_0^{N^{-\zeta}} \log |x|^{-1} d\nu_{B_N}^z(x) ]
&\le & \E[\log(\sigma^z_N)^{-1} 1_{\sigma^z_N\le N^{-\zeta}}]\\
&=&  C\left( N^{\frac{1}{2}+\kappa'-\zeta}
\right)^\beta \log (N^\zeta) +\int_0^{N^{-\zeta}} \frac{1}{x}
C\left( N^{\frac{1}{2}+\kappa'}
x\right)^\beta dx\\
\end{eqnarray*}
goes to zero as $N$ goes to infinity as soon
as $\zeta>\frac{1}{2}+\kappa'$.
We fix hereafter such a $\zeta$
and  we may and shall restrict the integration
from $N^{-\zeta}$ to $\e$.
To compare the integral for the spectral measure of $
A_N$ and $B_N$, observe that 
for all probability measure $P$, with $P_\gamma$ the
Cauchy law with parameter $\gamma$
\begin{equation}
	\label{ez1}
	P([a,b])\le P*P_\gamma([a-\eta,b+\eta])+P_\gamma([-\eta,\eta]^c)
\le  P*P_\gamma([a-\eta,b+\eta]) +\frac{\gamma}{\eta}
\end{equation}
whereas for $b-a>\eta$
\begin{equation}
	\label{ez2}
	P([a,b])\ge P*P_\gamma([a+\eta,b-\eta])-\frac{\gamma}{\eta}
	\,.
\end{equation}
Recall that
\begin{equation}
	\label{ez3}
	P*P_\gamma([a,b])= \int_a^b |\Im G(x+i\gamma)| dx.
\end{equation}
Set $\gamma=N^{-\kappa/5}$, $\kappa''=\kappa/2$ and $\eta=N^{-\kappa''/5}$.
We have, whenever $b-a\geq 4\eta$,
\begin{eqnarray*}
\E\nu^z_{B_N}([a,b])&\le& \int_{a-\eta}^{b+\eta} 
\E|\Im G_{B_n+z}(x+i\gamma)|
dx + N^{-(\kappa-\kappa'')/5}\\
&\le & (b-a+2N^{-\kappa''/5}
) + \nu^z_{A_N}*P_{N^{-\kappa/5}}( [a-N^{-\kappa/10},
b+N^{-\kappa/10}])+N^{-\kappa/10}\\
&\le& (b-a+2N^{-\kappa/10}) + \nu^z_{A_N}( [(a-2N^{-\kappa/10})_+,
(b+2N^{-\kappa/10})])+2N^{-\kappa/10}\,,
\end{eqnarray*}
where the first inequality is due to \eqref{ez1} and \eqref{ez3},
the second is due to 
Proposition \ref{prop1}, 
and the last uses \eqref{ez2} and \eqref{ez3}.
Therefore, if  $b-a=CN^{-\kappa/10}$ for some fixed $C$ larger than 4,
we deduce that there exists a finite constant $C'$ which only
depends on $C$ so that
$$\E\nu^z_{B_N}([a,b])\le C'(b-a)+ \nu^z_{A_N}( [(a-2N^{-\kappa/10})_+,
(b+2N^{-\kappa/10})])\,.$$
As a consequence, as we may assume without loss
of generality that $\kappa'>\kappa/10$,
\begin{eqnarray*}
&&\E[\int_{N^{-\zeta}}^\e\log |x|^{-1} d\nu_{B_N}^z(x)]\\
&\le&\sum_{k=0}^{[N^{\kappa/10}\e]}
\log(N^{-\zeta} +2Ck N^{-\kappa/10})^{-1} \E[\nu_{B_N}^z]([N^{-\zeta}
+2Ck N^{-\kappa/10}, N^{-\zeta}+2C(k+1) N^{-\kappa/10}])\,.
\end{eqnarray*}
We need to pay special attention to the
first term that we bound by noticing that
\begin{eqnarray*}
&&\log(N^{-\zeta} )^{-1} \E[\nu_{B_N}^z([N^{-\zeta}, 
N^{-\zeta}+2C N^{-\kappa/10}])]\\
&\le& \frac{10\zeta}{\kappa}
\log(N^{-\kappa/10} )^{-1} \E[\nu_{B_N}^z([0, 2(C+1) N^{-\kappa/10}])]\\
&\le& \frac{10\zeta}{\kappa}
\log(N^{-\kappa/10} )^{-1} (2C' N^{-\kappa/10}
+ \nu_{A_N}^z([0,(C+2) N^{-\kappa/10}]))\\
&\le&  \frac{20 C'\zeta}{\kappa} \log(N^{-\kappa/10} )^{-1}  N^{-\kappa/10}
+C''\int_0^{2(C+2) N^{-\kappa/10}}\log |x|^{-1} d \nu_{A_N}^z(x)\\
\end{eqnarray*}
For the other terms, we have
\begin{eqnarray*}
\sum_{k=1}^{[N^{\kappa/10}\e]}
&&\log(N^{-\zeta} +2Ck N^{-\kappa/10})^{-1} \E[\nu_{B_N}^z]([N^{-\zeta}
+2Ck N^{-\kappa/10}, N^{-\zeta}+2C(k+1) N^{-\kappa/10}])\\
&\le& 2C' \sum_{k=1}^{[N^{\kappa/10}\e]}
\log(Ck N^{-\kappa/10})^{-1} CN^{-\kappa/10}\\
&&+\sum_{k=1}^{[N^{\kappa/10}\e]}
\log(Ck N^{-\kappa/10})^{-1} \nu_{A_N}^z([2C(k-1) N^{-\kappa/10}, 
2C(k+2) N^{-\kappa/10}])\,.
\end{eqnarray*}
Finally, we can sum up all these inequalities
to find that there exists a finite constant
$C'''$ so that
$$\E[\int_{N^{-\zeta}}^\e\log |x|^{-1} d\nu_{B_N}^z(x)]
\le C'''\int_{0}^\e\log |x|^{-1} d\nu_{A_N}^z(x)
+C''' \int_0^\e \log |x|^{-1} dx$$
and therefore goes to zero when 
$n$ and then $\e$ goes to zero.
This proves the claim.
\qed
\section{Proof of Theorem \ref{theo2}.}
From the assumptions, it is clear that $(A_N+E_N)$ converges in $*$-moments
to the regular element $a$. By Theorem \ref{theo1}, it follows that
$L_N^{A+E+G}$
converges (weakly, in probability) towards $\nu_a$.
We can now remove $E_N$.
Indeed, by 
\eqref{eq-comp} and \eqref{eq-doublestar}, we have for any
$\chi<\kappa'/2$
and all $\xi\in\mathbb C$
$$|G^N_{A+G +\xi}(z)-G^N_{A+G+E+\xi}(z)|\le \frac{N^{-\chi}}{\Im z^2} $$
and therefore for $\Im z\ge N^{-\chi/2}$,
$$|\Im G^N_{A+G +\xi}(z)|\le |\Im G^N_{A+G+E+\xi}(z)| +1.$$
Again 
by \cite[Theorem 3.3]{SST} (or its generalization in \cite[Proposition 16]{GKZ})
to the complex case), 
for any $z$, the smallest singular value $\sigma^z_N$
of $A_N+G_N+z$ satisfies
$$P(\sigma^z_N\le x)\le C\left( N^{\frac{1}{2}+\kappa'}
x\right)^\beta$$
with $\beta=1$ or $2$ according whether we
are in the real or the complex case. We can now rerun the proof
of Theorem \ref{theo1}, replacing $A_N$ by $A_N'=A_N+E_N+G_N$ and
$B_N$ by $A_N'-E_N$.
\qed

\section{Proof of Corollary \ref{cor-tn}}
We apply Theorem \ref{theo2} with $A_N=T_N$, $E_N$ the $N$-by-$N$ matrix with
$$E_N(i,j)=\{\begin{array}{ll}
\delta_N=N^{-(1/2+\kappa')},& i=1,j=N\\
0,& \mbox{\rm otherwise}\,,
\end{array}
$$
where $\kappa'>\kappa$.
We check the assumptions of Theorem \ref{theo2}. We take $a$ to be a 
Unitary Haar element in ${\cal A}$, and recall that its Brown measure
$\nu_a$ is the uniform measure on $\{z\in\C: |z|=1\}$.
We now check that $a$ is regular.
Indeed, $\int x^k d\nu_a^z(x)=0$ if $k$ is odd by symmetry while
for $k$ even,
$$\int x^{k} d\nu_a^z(x)=\phi([(z-a)(z-a)^*]^{k/2})=
\sum_{j=1}^{k/2} (|z|^2+1)^{k-j}
\left(\begin{array}{l}
	k\\2j \end{array}\right)
\left(\begin{array}{l}
	2j\\j \end{array}\right)\,,
	$$
	and one therefore verifies that 
for $k$ even,
$$\int x^kd\nu_a^z(x)=\frac1{2\pi}\int (|z|^2+1+2|z|\cos \theta)^{k/2}
d\theta\,.$$ 
It follows that
$$\int_0^\e \log x d\nu_a^z(x)=\frac1{4\pi}\int_0^{2\pi} 
\log(|z|^2+1+2|z|\cos \theta){\bf 1}_{\{|z|^2+1+2|z|\cos\theta<\e\}}d\theta
\to_{\eps\to 0} 0\,,$$
proving the required regularity.

Further, we claim that 
$L_N^{A+E}$ converges to $\nu_a$.
Indeed the eigenvalues $\lambda$ of $A_N+E_N$
are such that there exists a non-vanishing vector
$u$ so that 
$$u_N\delta_N=\lambda  u_1,u_{i-1}=\lambda u_i\,,$$
that is
$$\lambda^N=\delta_N.$$
In particular, all the $N$-roots of $\delta_N$ are (distinct) eigenvalues, that
is the eigenvalues $\lambda_j^N$ of $A_N$ are
$$\lambda_j^N=|\delta_N|^{1/N} e^{2i\pi j/N},\quad 1\le j\le N\,.$$
Therefore,
 for any bounded continuous $g$ function on $\C$,
$$ \lim_{N\ra\infty} \frac1N \sum_{i=1}^N  g(\lambda_j^N)=\frac1{2\pi}
\int g(\theta) d\theta\,,$$
as claimed.
\qed

\section{Proof of Proposition~\ref{prop-lb}}
 \label{sec-proppf}

In this section we will prove the following proposition:

\begin{proposition}\label{prop-lb2}
  Let $b= b(N)$ be a sequence of positive integers, and let $T_{b,N}$ be as in Proposition~\ref{prop-lb}.  Let $R_N$ be an $N$ by $N$ matrix satisfying $\| R_N\| \le g(N)$, where for all $N$ we assume that $g(N) < \frac1{3 b \sqrt N}$.  Then
$$ \rho( T_{b,N} + R_N) \le \left( O\left( \sqrt{Nb} \left(2 N^{1/4} g^{1/2}\right)^b\right)\right)^{1/(b+1)} + \left( b^2 Ng\right)^{1/(b+1)}\,.$$
\end{proposition}

Proposition~\ref{prop-lb} follows from Proposition~\ref{prop-lb2} by adding the assumption that $b(N) \ge log(N)$ and then simplifying the upper bound on the spectral radius.

\bigskip
\noindent
{\bf Proof of Proposition~\ref{prop-lb2}:}
To bound the spectral radius, we will use the fact that $\rho( T_{b,N} + R_N) \le \left\| (T_{b,N} + R_N)^k \right\|^{1/k}$ for all integers $k \ge 1$.  Our general plan will be to bound $\left\| (T_{b,N} + R_N)^k \right\|$ and then take a $k$-th root of the bound.  We will take $k=b+1$, which allows us to take advantage of the fact that $T_{b,N}$ is $(b+1)$-step nilpotent.  In particular, we make use of the fact that for any positive integer $a$,
\begin{equation}
 \label{eqn-tbna}
\| \tbn^a \| = \begin{cases}
                (b-a+1)^{1/2} \floor{\frac{N}{b+1}}^{1/2} & \mbox{ if } 1 \le a \le b\\
		0 & \mbox{ if } b+1 \le a.
               \end{cases}
\end{equation}

We may write
\begin{align}
\norm{ (\tbn + R_N)^{b+1} } 
&\le \sum_{\lambda \in \{0,1\}^{b+1}} \norm{ \prod_{i=1}^{b+1} \tbn^{\lambda_i} R_N^{1-\lambda_i} }\nonumber\\
&= \sum_{\ell =0}^{b+1} \mathop{\sum_{\lambda \in \{0,1\}^{b+1}}}_{\lambda \mathrm{\ has\ } \ell \mathrm{\ ones}} \norm{ \prod_{i=1}^{b+1} \tbn^{\lambda_i} R_N^{1-\lambda_i} }\nonumber
\end{align}
When $\ell$ is large, we will make use of the following lemma.

\begin{lemma}
 \label{lem-secsum1}
If $\lambda\in \{0,1\}^{k}$ has $\ell$ ones and $\ell \ge (k+1)/2$, then
$$\norm{ \prod_{i=1}^{k} \tbn^{\lambda_i} R_N^{1-\lambda_i} } 
\le \norm{ \tbn^{\floor{\frac{\ell}{k-\ell+1}}}}^{k-\ell+1} \norm{R_N}^{k-\ell}.
$$
\end{lemma}

We will prove Lemma~\ref{lem-secsum1} in Section~\ref{lempfs}.

Using Lemma~\ref{lem-secsum1} with $k=b+1$ along with the fact that $\norm{AB}\le \norm A \norm B$, we have
\begin{align}
\norm{ (\tbn + R_N)^{b+1} } 
&\le \sum_{\ell =0}^{\floor{\frac{b+2}{2}}} \binom{b+1}{\ell} \norm{\tbn}^\ell \norm{R_n}^{b-\ell+1} \nonumber\\
&\qquad
   + \sum_{\ell=\ceiling{\frac{b+2}{2}}}^{b+1} \binom{b+1}{\ell} \norm{ \tbn^{\floor{\frac{\ell}{b-\ell+2}}}}^{b-\ell+2} \norm{R_N}^{b-\ell+1}.\nonumber\\
&\le \sum_{\ell =0}^{\floor{\frac{b+2}{2}}} \binom{b+1}{\ell} \norm{\tbn}^\ell g^{b-\ell+1} \label{firstsum}\\
&\qquad
   + \sum_{\ell=\ceiling{\frac{b+2}{2}}}^{b+1} \binom{b+1}{\ell} \norm{ \tbn^{\floor{\frac{\ell}{b-\ell+2}}}}^{b-\ell+2} g^{b-\ell+1},\label{secondsum}
\end{align}
where the second inequality comes from the assumption $\|R_N\|\le g = g(N)$.

We will bound \eqref{firstsum} and \eqref{secondsum} separately.  To bound \eqref{firstsum} note that
\begin{align}
 \sum_{\ell =0}^{\floor{\frac{b+2}{2}}} \binom{b+1}{\ell} \norm{\tbn}^\ell g^{b-\ell+1}
&\le 
\sum_{\ell =0}^{\floor{\frac{b+2}{2}}} \binom{b+1}{\ell} \left( (b+1)\floor{\frac{N}{b+1}}\right)^{\ell/2} g^{b-\ell+1} \nonumber\\
&\le \frac{b+4}{2} \binom{b+1}{\floor{(b+1)/2}} N^{(b+2)/4} g^{b/2} \nonumber\\
&=O\left(\sqrt{Nb}(2N^{1/4}g^{1/2})^b\right).\label{fsbd}
\end{align}

Next, we turn to bounding \eqref{secondsum}. We will use the following lemma to show that the largest term in the sum \eqref{secondsum} comes from the $\ell=b$ term.  Note that when $\ell=b+1$, the summand in \eqref{secondsum} is equal to zero by \eqref{eqn-tbna}.

\begin{lemma};
 \label{lem-secsum}
  If $\norm{ \tbn^{\floor{\frac{\ell+1}{b-\ell+1}}}} > 0$ and $\ell \le b-1$ 
and $$ g \le \frac{2}{e^{3/2} N^{1/2}b}, $$
then
$$
\binom{b+1}{\ell} \norm{ \tbn^{\floor{\frac{\ell}{b-\ell+2}}}}^{b-\ell+2} g^{b-\ell+1}
\le
\binom{b+1}{\ell+1} \norm{ \tbn^{\floor{\frac{\ell+1}{b-\ell+1}}}}^{b-\ell+1} g^{b-\ell}.
$$
\end{lemma}

We will prove Lemma~\ref{lem-secsum} in Section~\ref{lempfs}.

Using Lemma~\ref{lem-secsum} we have
\begin{align}
 \sum_{\ell=\ceiling{\frac{b+2}{2}}}^{b+1} \binom{b+1}{\ell} \norm{ \tbn^{\floor{\frac{\ell}{b-\ell+2}}}}^{b-\ell+2} g^{b-\ell+1} 
&\le \frac b2 (b+1) \norm{ \tbn^{\floor{\frac{b}{2}}}}^{2} g^{1} \nonumber\\
&\le \frac b2 (b+1) (b-\floor{b/2} +1)\frac{N}{b+1} g \nonumber\\
&\le b^2 Ng. \label{ssbd}
\end{align}

Combining \eqref{fsbd} and \eqref{ssbd} with \eqref{firstsum} and \eqref{secondsum}, we may use the fact that $(x+y)^{1/(b+1)} \le x^{1/(b+1)} +y^{1/(b+1)}$ for positive $x,y$ to complete the proof of Proposition~\ref{prop-lb2}.  It remains to prove Lemma~\ref{lem-secsum1} and Lemma~\ref{lem-secsum}, which we do in Section~\ref{lempfs} below.
\qed

\subsection{Proofs of Lemma~\ref{lem-secsum1} and Lemma~\ref{lem-secsum}}\label{lempfs}

\bigskip
\noindent
{\bf Proof of Lemma~\ref{lem-secsum1}:}
Using \eqref{eqn-tbna}, it is easy to show that 
\begin{equation}
 \label{relation1}
  \norm{\tbn^a}\norm{\tbn^c} < \norm{\tbn^{a-1}}\norm{\tbn^{c+1}} \mbox{ for integers } 3\le c+2 \le a \le b.
\end{equation}
It is also clear from \eqref{eqn-tbna} that 
\begin{equation}
 \label{relation2}
\norm{\tbn^a} \le \norm{\tbn^{a-1}} \mbox{ for all positive integers $a$}.
\end{equation}

Let $\lambda \in \{0,1\}^{k}$ have $\ell$ ones.  Then, using the assumption that $\ell \ge k-\ell +1$, we may write
$$\prod_{i=1}^{k} \tbn^{\lambda_i} R_N^{1-\lambda_i} 
= \tbn^{a_1} R_N^{b_1}
\tbn^{a_2} R_N^{b_2}
\cdots\tbn^{a_{k-\ell}} R_N^{b_{k-\ell}}\tbn^{a_{k-\ell+1}},$$
where $a_i \ge 1$ for all $i$ and $b_i \ge 0$ for all $i$.
Thus
$$\norm{\prod_{i=1}^{k} \tbn^{\lambda_i} R_N^{1-\lambda_i}} \le
\norm{R_N}^{k-\ell} \prod_{i=1}^{k-\ell+1} \norm{\tbn^{a_i}}.$$
Applying \eqref{relation1} repeatedly, we may assume that two of the $a_i$ differ by more than 1, all without changing the fact that $\sum_{i=1}^{k-\ell+1}a_i = \ell$.  Thus, some of the $a_i$ are equal to $\floor{\frac{\ell}{k-\ell+1}}$ and some are equal to $\ceiling{\frac{\ell}{k-\ell+1}}$.  Finally, applying \eqref{relation2}, we have that
$$\prod_{i=1}^{k-\ell+1} \norm{\tbn^{a_i}} \le \norm{\tbn^{\floor{\frac{\ell}{k-\ell+1}}}}^{k-\ell+1}.$$
\qed

\bigskip
\noindent
{\bf Proof of Lemma~\ref{lem-secsum}:}
Using \eqref{eqn-tbna} and rearranging, it is sufficient to show that
$$\frac{\ell+1}{b-\ell+1} \left(b -\floor{\frac{\ell}{b-\ell+2}} +1\right) ^{1/2} \floor{\frac{N}{b+1}}^{1/2} g 
\le \left(  
\frac{b-\floor{\frac{\ell+1}{b-\ell+1}}+1}
{b-\floor{\frac{\ell}{b-\ell+2}}+1}
\right)^{\frac{b-\ell+1}{2}}
$$
Using a variety of manipulations, it is possible to show that 
\begin{eqnarray*}
\left(  
\frac{b-\floor{\frac{\ell+1}{b-\ell+1}}+1}
{b-\floor{\frac{\ell}{b-\ell+2}}+1}
\right)^{\frac{b-\ell+1}{2}}
&\ge& 
\exp\left(
-\frac{(b-\ell+2)(b-\ell+1)}
{(b+2)(b-\ell+2)-\ell}
-
\frac{b+2}
{(b+2)(b-\ell+2)-\ell}
\right)\\&
\ge&\exp(-3/2).
\end{eqnarray*}
Thus, it is sufficient to have 
$$\frac b2 N^{1/2} g \le \exp(-3/2),$$
which is true by assumption.
\qed

\bibliographystyle{amsplain}

\end{document}